\newif\ifpdf
\ifx\pdfoutput\undefined
\pdffalse 
\else
\pdfoutput=1 
\pdftrue \fi

\newif\iffinal
\finalfalse	
\finaltrue 

\documentclass[reqno,twoside,11pt]{amsart}
\iffinal\else\usepackage[notref,notcite]{showkeys}\fi
\usepackage{amsmath}
\usepackage{amsfonts}
\usepackage{amssymb}
\usepackage{verbatim}
\usepackage{epsfig}

\IfFileExists{epsf.def}{\input epsf.def}{\usepackage{epsf}}

\usepackage{mathpazo}
\usepackage{mathrsfs}

\DeclareFontFamily{OT1}{eusb}{} \DeclareFontShape{OT1}{eusb}{m}{n} {<5> <6> <7> <8> <9> <10> <11> <12> <14.4> eusb10}{}
\DeclareMathAlphabet{\eusb}{OT1}{eusb}{m}{n}

\DeclareFontFamily{OT1}{eusm}{} \DeclareFontShape{OT1}{eusm}{m}{n} {<5> <6> <7> <8> <9> <10> <11> <12> <14.4> eusm10}{}
\DeclareMathAlphabet{\eusm}{OT1}{eusm}{m}{n}

\DeclareFontFamily{OT1}{eufm}{} \DeclareFontShape{OT1}{eufm}{m}{n} {<5> <6> <7> <8> <9> <10> <11> <12> <14.4> eufm10}{}
\DeclareMathAlphabet{\mathfrak}{OT1}{eufm}{m}{n}

\DeclareFontFamily{OT1}{fraktura}{}
\DeclareFontShape{OT1}{fraktura}{m}{n} {<5> <6> <7> <8> <9> <10> <11> <12> <13> <14.4> [1.1] eufm10}{}
\DeclareMathAlphabet{\fraktura}{OT1}{fraktura}{m}{n}

\DeclareFontFamily{OT1}{cmfi}{} \DeclareFontShape{OT1}{cmfi}{m}{n} {<5> <6> <7> <8> <9> <10> <11> <12> <13> <14.4> [0.9] cmfi10}{}
\DeclareMathAlphabet{\cmfi}{OT1}{cmfi}{b}{n}

\DeclareFontFamily{OT1}{cmss}{} \DeclareFontShape{OT1}{cmss}{m}{n} {<5> <6> <7> <8> <9> <10> <11> <12> <13> <14.4> cmss10}{}
\DeclareMathAlphabet{\cmss}{OT1}{cmss}{m}{n}

\newtheoremstyle{thm}{1.5ex}{1.5ex}{\itshape\rmfamily}{} {\bfseries\rmfamily}{}{2ex}{}

\newtheoremstyle{def}{1.5ex}{1.5ex}{\slshape\rmfamily}{} {\bfseries\rmfamily}{}{2ex}{}

\newtheoremstyle{rem}{1.5ex}{1.5ex}{\rmfamily}{} {\bfseries\rmfamily}{}{2ex}{}

\newenvironment{proofsect}[1] {\vskip0.1cm\noindent{\rmfamily\itshape#1.}}{\qed\vspace{0.15cm}}

\theoremstyle{thm}
\newtheorem{theorem}{Theorem}[section]
\newtheorem{lemma}[theorem]{Lemma}
\newtheorem{proposition}[theorem]{Proposition}

\newtheorem*{Main Theorem}{Main Theorem.}
\newtheorem{corollary}[theorem]{Corollary}
\newtheorem*{special theorem}{Lindeberg-Feller Theorem for Martingales}

\theoremstyle{def}
\newtheorem{definition}[theorem]{Definition}

\theoremstyle{rem}
\newtheorem{remark}[theorem]{Remark}

\numberwithin{equation}{section}

\renewcommand{\section}{\secdef\sct\sect}
\newcommand{\sct}[2][default]{%
\refstepcounter{section}
\addcontentsline{toc}{section}{{\tocsection {}{\thesection}{\!\!\!\!#1\dotfill}}{}}
\vspace{0.7cm}
\centerline{\scshape\thesection.\ #1} \nopagebreak \vspace{0.2cm}}
\newcommand{\sect}[1]{%
\vspace{0.4cm} \centerline{\large\scshape\rmfamily #1}
\vspace{0.2cm}}

\renewcommand{\subsection}{\secdef\subsct\sbsect}
\newcommand{\subsct}[2][default]{\refstepcounter{subsection}
\addcontentsline{toc}{subsection}
{{\tocsection{\!\!}{\hspace{1.2em}\thesubsection}{\!\!\!\!#1\dotfill}}{}}
\nopagebreak\vspace{0.45\baselineskip} {\flushleft\bf
\thesubsection~\bf #1.~}
\\*[3mm]\noindent
\nopagebreak}
\newcommand{\sbsect}[1]{\vspace{0.1cm}\noindent
\textbf{#1.~}\vspace{0.1cm}}

\renewcommand{\subsubsection}{%
\secdef \subsubsect\sbsbsect}
\newcommand{\subsubsect}[2][default]{%
\refstepcounter{subsubsection} 
\addcontentsline{toc}{subsubsection}{{\tocsection{\!\!}
{\hspace{3.05em}\thesubsubsection}{\!\!\!\!#1\dotfill}}{}}
\nopagebreak
\vspace{0.15\baselineskip} \nopagebreak {\flushleft\rmfamily
\itshape\thesubsubsection
\ \rmfamily #1\/.}\ }
\newcommand{\sbsbsect}[1]{\vspace{0.1cm}\noindent
\rmfamily \itshape
\arabic{section}.\arabic{subsection}.\arabic{subsubsection} \
\sffamily #1\/.\ }

\iffinal

\else

\fi

\renewcommand{\caption}[1]{%
\vglue0.5cm
\refstepcounter{figure}
\begin{minipage}{0.9\textwidth}\small {\sc Figure~\thefigure. }#1\end{minipage}}

\newcommand{\diam}{\operatorname{diam}}

\newcommand{\texte}{\text{\rm e}}

\renewcommand{\AA}{\mathcal A}

\newcommand{\EE}{\mathcal E}
\newcommand{\FF}{\mathcal F}
\newcommand{\GG}{\mathcal G}

\newcommand{\E}{\mathbb E}

\newcommand{\G}{\mathbb G}

\newcommand{\BbbP}{\mathbb P}

\newcommand{\R}{\mathbb R}

\newcommand{\Z}{\mathbb Z}

\newcommand{\scrB}{\mathscr{B}}

\newcommand{\twoeqref}[2]{(\ref{#1}--\ref{#2})}
\newcommand{\cc}{{\text{\rm c}}}

\def\myffrac#1#2 in #3{\raise 2.6pt\hbox{$#3 #1$}\mkern-1.5mu\raise 0.8pt\hbox{$#3/$}\mkern-1.1mu\lower 1.5pt\hbox{$#3 #2$}}
\newcommand{\ffrac}[2]{\mathchoice%
	{\myffrac{#1}{#2} in \scriptstyle}
	{\myffrac{#1}{#2} in \scriptstyle}
	{\myffrac{#1}{#2} in \scriptscriptstyle}
	{\myffrac{#1}{#2} in \scriptscriptstyle}
}

\begin{document}

\title[Graph diameter in long-range percolation]
{Graph diameter in long-range percolation}

\author[Marek~Biskup]
{Marek~Biskup}

\thanks{\hglue-4.5mm\fontsize{9.6}{9.6}\selectfont\copyright\,2009 by M.~Biskup. Reproduction, by any means, of the entire article for non-commercial purposes is permitted without charge.\vspace{2mm}}
\maketitle

\vspace{-5mm}
\centerline{\textit{Department of Mathematics, University of California at Los Angeles, USA}}
\centerline{\textit{School of Economics, University of South Bohemia, Czech Republic}}

\vspace{2mm}
\begin{quote}
\footnotesize \textbf{Abstract:}
We study the asymptotic growth of the diameter of a graph obtained by adding sparse ``long'' edges to a square box in~$\Z^d$. We focus on the cases when an edge between~$x$ and~$y$ is added with probability decaying with the Euclidean distance as~$|x-y|^{-s+o(1)}$ when~$|x-y|\to\infty$. For~$s\in(d,2d)$ we show that the graph diameter for the graph reduced to a box of side~$L$ scales like~$(\log L)^{\Delta+o(1)}$ where~$\Delta^{-1}:=\log_2(2d/s)$. In particular, the diameter grows about as fast as the typical graph distance between two vertices at distance~$L$. We also show that a ball of radius~$r$ in the intrinsic metric on the (infinite) graph will roughly coincide with a ball of radius $\exp\{r^{1/\Delta+o(1)}\}$ in the Euclidean metric.
\end{quote}
\vspace{2mm}

\section{Main result}
\noindent
Consider the~$d$-dimensional hypercubic lattice~$\Z^d$ and add a random collection of edges to~$\Z^d$ according to the following rule: An edge between distinct sites~$x$ and~$y$ occurs with probability~$p_{xy}$, independently of all other edges, where~$p_{xy}$ depends only on the difference~$x-y$ and decays like~$|x-y|^{-s+o(1)}$ as the Euclidean norm $|x-y|$ tends to infinity. Let~$D(x,y)$ denote the graph distance between~$x$ and~$y$ which is defined as the length of the shortest path that connects~$x$ to~$y$ using only edges that are available in the present (random) sample.

In \cite{Biskup} we studied the asymptotic of~$D(x,y)$ as~$|x-y|\to\infty$. In particular, it was shown that for~$s\in(d,2d)$ this distance behaves like
\begin{equation}
\label{Dxy} 
D(x,y)=\bigl(\log|x-y|\bigr)^{\Delta+o(1)},
\qquad|x-y|\to\infty,
\end{equation}
where
\begin{equation}
\label{Delta} 
\Delta:=\frac {\log2}{\log\bigl(\frac{2d}s\bigr)}.
\end{equation}
Technically, \eqref{Dxy} is established with ``$o(1)$ tending to zero in probability'' and thus represents the \emph{typical} behavior for fixed~$x$ and~$y$. The result allows for the possibility that even the nearest-neighbor edges are randomized --- $x$ and~$y$ are then restricted to the unique infinite connected component.

The main purpose of this note is to determine the corresponding asymptotic for the \emph{maximal} graph distance between any two sites in a large, finite set. Explicitly, let us consider the box~$\Lambda_L:=[-L,L]^d\cap\Z^d$ and let $\G_L$ denote the restriction of the above random graph to vertices, and edges with both endpoints, in~$\Lambda_L$. Let~$D_L(x,y)$ denote graph-theoretical (a.k.a.~intrinsic or chemical) distance between~$x,y\in\Lambda_L$ as measured on~$\G_L$. The diameter of~$\G_L$ is then given by
\begin{equation}
D_L:=\max_{x,y\in\Lambda_L}D_L(x,y).
\end{equation}
The following settles a question that was left open in~\cite{Biskup}:

\begin{theorem}
\label{thm1.1}
Suppose that~$p_{xy}$ can be written as~$p_{xy}:=1-\texte^{-q(x-y)}$, where~$q\colon\Z^d\to[0,\infty]$ is an even function for which the limit
\begin{equation}
\label{plim}
s:=-\lim_{|x|\to\infty}\,\frac{\log q(x)}{\log|x|}
\end{equation}
exists and satisfies~$s\in(d,2d)$. Then for all~$\epsilon>0$,
\begin{equation}
\label{lim}
\lim_{L\to\infty}\BbbP\bigl((\log L)^{\Delta-\epsilon}\le D_L\le(\log L)^{\Delta+\epsilon}\bigr)=1,
\end{equation}
where~$\Delta$ is as in \eqref{Delta}.
\end{theorem}

It is clear that the asymptotic in \eqref{Dxy} serves as a lower bound on~$D_L$. However, a matching upper bound --- the main contribution of this note --- turns out to be much less immediate. The point is that the bounds from \cite{Biskup} for $\BbbP(D(x,y)\ge(\log|x-y|)^{\Delta+\epsilon})$ are much too weak to imply the same upper bound on~$D_L$. This is because the strategy employed in~\cite{Biskup} is based on the fact that one can find, with overwhelming probability, long edges within up to $(\log L)^{o(1)}$-distance from any given site. Unfortunately, this does not hold for every site of~$\Lambda_L$; in fact, $\Lambda_L$ \emph{will} contain a translate of~$\Lambda_\ell$ with~$\ell\approx(\log L)^{1/d}$ whose vertices have \emph{no other} edges than those inherited from~$\Z^d$.

The restriction of the above result to hypercubic lattice is mostly a matter of convenience; no part of the proof depends essentially on the details of the underlying graph. (What we need is that the graph is embedded in~$\R^d$ so that \emph{any} square block of side~$L$ contains order~$L^d$ sites.) Similarly, we could also work with more general sequences of sets than cubic boxes. In fact, we could even accommodate non-translation invariant distributions and/or diluted lattices, e.g., work (as was done in~\cite{Biskup}) with long-range percolation under the sole assumption that there is an infinite connected component. Notwithstanding, such generalizations tend to obscure the main ideas of the proof and so we settle to a translation-invariant and \emph{a priori} connected setting.

\smallskip
The control of the diameter provided by Theorem~\ref{thm1.1} allows for some control of the volume growth of the percolation graph. Consider a realization of the long-range percolation on~$\Z^d$ --- still including the edges of~$\Z^d$ --- and let
\begin{equation}
B(0,r):=\bigl\{x\in\Z^d\colon D(0,x)\le r\bigr\}
\end{equation}
denote the ball of radius~$r$ in the intrinsic metric. Trapman~\cite[Theorem~1.1(c)]{Trapman} has recently shown that the volume of this ball grows subexponentially with the radius, i.e.,
\begin{equation}
s>d\quad\Rightarrow\quad \lim_{r\to\infty}\bigl|B(0,r)|^{1/r}=1,\quad\BbbP\text{\rm-a.s.}
\end{equation}
Here we derive the leading order of the growth of $|B(0,r)|$ with~$r$:

\begin{theorem}
\label{thm1.2}
Under the conditions of Theorem~\ref{thm1.1}, for each~$\epsilon>0$,
\begin{equation}
\label{eq1.8}
\lim_{r\to\infty}\BbbP\bigl(\Lambda_{\,\exp\{r^{-\epsilon+1/\Delta}\}}\subset B(0, r)\subset\Lambda_{\,\exp\{r^{\,\epsilon+1/\Delta}\}}\bigr)=1.
\end{equation}
In particular,
\begin{equation}
\frac{\log\log |B(0,r)|}{\log r}\,\underset{r\to\infty}\longrightarrow\,\frac1\Delta
\end{equation}
in probability.
\end{theorem}

As $\Delta\in(1,\infty)$ for $s\in(d,2d)$, the leading-order volume growth takes a stretched-exponential form, i.e.,
\begin{equation}
\bigl|B(0,r)\bigr|=\exp\bigl\{ r^{\frac1\Delta+o(1)}\bigr\},\qquad r\to\infty.
\end{equation}
While the left inclusion in \eqref{eq1.8} is implied directly by Theorem~\ref{thm1.1}, for that on the right we will have to invoke --- and, in fact, prove again in order to accommodate for a more general setting --- a result due to Trapman (see Theorem~\ref{thm3.1}).

\smallskip
The rest of this note is organized as follows: In Sect.~\ref{sec2} we discuss various motivations for, and further results related to this work. In Sect.~\ref{sec3} we prove Theorem~\ref{thm1.2} concerning volume growth estimates on the infinite graph. Sect.~\ref{sec4} gives the proof of Theorem~\ref{thm1.1} on graph diameter subject to some technical claims; these are then established in Sect.~\ref{sec5}.

\section{Related work}
\label{sec2}\noindent
Long-range percolation, of which our model is an example, originated in the mathema\-ti\-cal-physics literature as a model that exhibits a phase transition even in spatial dimension one (e.g., Newman and Schulman~\cite{Newman-Schulman}, Schulman~\cite{Schulman}, Aizenman and Newman~\cite{Aizenman-Newman}, Imbrie and Newman~\cite{Imbrie-Newman}). It soon became clear that $s=d$ and $s=2d$ are two distinguished values; for $s<d$ the model is essentially mean-field (or complete-graph) alike, for $s>2d$ the behavior is more or less as for the nearest-neighbor percolation. The regime $d<s<2d$ turned out to be quite interesting; indeed, it is the only general class of percolation models with Euclidean (or amenable) geometry where one can prove absence of percolation at the percolation threshold (Berger~\cite{Berger-RW}).
In all dimensions, the model with $s=2d$ has a natural continuum scaling limit.

Recently, long-range percolation has been invoked as a fruitful source of graphs with non-trivial growth properties. Our interest was stirred by the work of Benjamini and Berger~\cite{Benjamini-Berger} who proposed (and studied) long-range percolation as a model of social networks. It is  this context where the graph distance scaling, and volume growth, are particularly of much interest. Thanks to numerous contributions that followed \cite{Benjamini-Berger}, this scaling is now known for most values of~$s$ and~$d$. Explicitly, for~$s<d$, a corollary to the main result of Benjamini, Kesten, Peres and Schramm \cite{BKPS} asserts that
\begin{equation}
D_L\,\underset{L\to\infty}\longrightarrow\,\Bigl\lceil\frac s{d-s}\Bigr\rceil,
\end{equation}
almost surely.
As~$s\uparrow d$, the right-hand side tends to infinity and so, at~$s=d$, we expect~$D_L\to\infty$. And, indeed, the precise growth rate in this case has been established by Coppersmith, Gamarnik and Sviridenko \cite{CGS},
\begin{equation}
\label{DLs=d}
D_L\asymp\frac{\log L}{\log\log L},\qquad L\to\infty,
\end{equation}
where~``$\asymp$'' means that the ratio of left and right-hand side is a random variable that is bounded away from zero and infinity with probability tending to one.

For $s\in(d,2d)$, the present paper states $D_L=(\log L)^{\Delta+o(1)}$. Here we note that~$\Delta\downarrow1$ as~$s\downarrow d$ which, formally, is in agreement with \eqref{DLs=d}. For~$s\uparrow2d$ we in turn have~$\Delta\to\infty$ and so, at~$s=2d$, a polylogarithmic growth is no longer sustainable. Instead, for the case of the decay~$p_{xy}\sim\beta|x-y|^{-2d}$ one expects that
\begin{equation}
\label{s=2d}
D_L=L^{\theta(\beta)+o(1)},\qquad L\to\infty,
\end{equation}
where~$\theta(\beta)$ varies through~$(0,1)$ as $\beta$ sweeps through $(0,\infty)$. This claim is supported by upper and lower bounds in somewhat restricted one-dimensional cases (Benjamini and Berger \cite{Benjamini-Berger}, Coppersmith, Gamarnik and Sviridenko \cite{CGS}). However, even the existence of a sharp exponent $\theta(\beta)$ has been elusive so far.

For~$s>2d$ one expects \cite{Benjamini-Berger} the same behavior as for the original graph. And indeed, the linear asymptotic,
\begin{equation}
D_L\asymp L,
\end{equation}
 has been established by Berger \cite{Berger-LRP}. For the nearest-neighbor percolation case, this statement goes back to the work of Antal and Pisztora \cite{Antal-Pisztora}.

Further motivation comes from the recent interest in diffusive properties of graphs arising via long-range percolation. An early work in this respect was that of Berger~\cite{Berger-RW} who characterized regimes of recurrence and transience for the simple random walk on such graphs. Benjamini, Berger and Yadin~\cite{Benjamini-Berger-Yadin} later showed that the mixing time $\tau_L$ of the random walk on~$\G_L$ in $d=1$ scales like 
\begin{equation}
\tau_L\sim
\begin{cases}
L^{s-1},\qquad&\text{if }1<s<2,
\\
L^2,\qquad&\text{if }s=2,
\end{cases}
\end{equation}
with an apparent jump in the exponent when~$s$ passes through~2.  Misumi~\cite{Misumi} found estimates on the effective resistance in $\Lambda_{2L}\setminus\Lambda_L$ that exhibit a similar transition. 

Very recently, precise bounds for the heat kernel and spectral gap of such random walks have been derived by Crawford and Sly~\cite{Crawford-Sly}. These are claimed to lead to the proof that the law of such random walks scales to $\alpha$-stable processes for $d<s<d+2$ in $d\ge2$ and $1<s<2$ in $d=1$. For $s$ on the increasing side of these regimes, the random walk is expected to scale to Brownian motion.

\section{Volume growth}
\label{sec3}\noindent
The goal of this section is to prove Theorem~\ref{thm1.2}. As already mentioned, while the left inclusion in \eqref{eq1.8} is a direct consequence of Theorem~\ref{thm1.1}, the proof of the right inclusion will be based on ideas underlying the proof of Theorem~1.2 in Trapman~\cite{Trapman}. Unfortunately, Trapman's setting is too stringent for our purposes and so we restate (and prove) the relevant result in a more suitable form:

\begin{theorem}
\label{thm3.1}
Under the conditions of Theorem~\ref{thm1.1}, for each $s'\in(d,s)$ there are constants $c_1,c_2\in(0,\infty)$ such that, for $\Delta':=1/\log_2(2d/s')$,
\begin{equation}
\label{my-bd}
\BbbP\bigl(D(0,x)\le n\bigr)\le c_1\biggl(\frac{\texte^{\,c_2\, n^{1/\Delta'}}}{|x|}\biggr)^{s'},\qquad n\ge1.
\end{equation}
\end{theorem}

Before we provide a proof of this result, let us see how it fits into our proof of the volume growth estimate:

\begin{proofsect}{Proof of Theorem~\ref{thm1.2}}
Notice first that, by the structure of the expressions, it suffices to prove both limits in the statement along a single sequence of~$r$'s that tends to infinity at most exponentially fast. In fact, we will do this for $r$ being of the form $(\log L)^\theta$ with~$\theta\approx\Delta$ and~$L$ is running through positive integers.

We begin with the right inclusion in~\eqref{eq1.8}. Let $\epsilon>0$ and pick $s'\in(d,s)$ so that $\Delta':=1/\log_2(2d/s')$ satisfies $\Delta'>\Delta-\epsilon$. Setting $\beta:=(\Delta-\epsilon)/\Delta'$, a union bound and Theorem~\ref{thm3.1} then give
\begin{multline}
\qquad
\BbbP\bigl(\exists x\in\Lambda_{2^{k+1}L}\setminus\Lambda_{2^k L}\colon D(0,x)\le(\log L)^{\Delta-\epsilon}\bigr)
\\
\le
c_3\biggl(\frac{\exp\{c_2(\log L)^\beta\}}{2^kL}\biggr)^{s'}\bigl(2^{k+1} L\bigr)^d
=L^{d-s'+o(1)} (2^k)^{d-s'},
\qquad
\end{multline}
where we used that $\beta<1$ by our assumptions and where $o(1)\to0$ in the limit as $L\to\infty$. Since~$s'>d$, the right-hand side is summable on~$k$ and so we conclude
\begin{equation}
\label{eq3.3}
\BbbP\bigl(\exists x\not\in\Lambda_L\colon D(0,x)\le(\log L)^{\Delta-\epsilon}\bigr)\le L^{d-s'+o(1)},
\end{equation}
which tends to zero as $L\to\infty$. A moment's thought shows that
\begin{equation}
\bigl\{\exists x\not\in\Lambda_L\colon D(0,x)\le(\log L)^{\Delta-\epsilon}\bigr\}
\supset\bigl\{B(0,r)\not\subset\Lambda_{\,\exp\{r^{\,\epsilon'+1/\Delta}\}\,}\bigr\}
\end{equation}
for $r:=(\log L)^{\Delta-\epsilon}$ and $\epsilon':=(\Delta-\epsilon)^{-1}-\Delta^{-1}$. The right inclusion in \eqref{eq1.8} thus holds with probability tending to one for all $\epsilon>0$.

As to the left inclusion in \eqref{eq1.8} we notice that for $r:=(\log L)^{\Delta+\epsilon}$,
\begin{equation}
\bigl\{D_L\le(\log L)^{\Delta+\epsilon}\bigr\}\subset \bigl\{\Lambda_{\,\exp\{r^{-\epsilon'+1/\Delta}\}}\subset B(0,r)\bigr\}
\end{equation}
where $\epsilon':=\Delta^{-1}-(\Delta+\epsilon)^{-1}$. By Theorem~\ref{thm1.1}, the event on the left occurs with probability tending to one as~$L\to\infty$. Therefore, so does the left inclusion in~\eqref{eq1.8}.
\end{proofsect}

In order to prove Theorem~\ref{thm3.1}, we will follow Trapman's remarkable simplification of the proof from Biskup~\cite{Biskup} for the lower bound on the graph distance in infinite-volume setting. Fix $s'\in(d,s)$ and let~$R=R(s')\ge1$ be the number such that
\begin{equation}
\label{2.3a}
p_{xy}\le |x-y|^{-s'},\qquad |x-y|\ge R.
\end{equation}
This number exists by our assumption \eqref{plim}. The key steps of Trapman's argument can be encapsulated into two lemmas:

\begin{lemma}
\label{lemma-step1}
Abbreviate $B_k:=B(0,k)$. If $|x|/k\ge R$, then
\begin{equation}
\label{PB-bd}
\BbbP\bigl(D(0,x)\le k\bigr)\le\Bigl(\frac{|x|}k\Bigr)^{-s'}\sum_{j=0}^k \E|B_{j}|\,\E|B_{n-j}|
\end{equation}
\end{lemma}

\begin{proofsect}{Proof}
If $D(0,x)\le k$, then there exists a (vertex) self-avoiding path from $0$ to~$x$ such that at least one edge has length at least $|x|/k$. If this edge occurs at the~$j$-th step and it goes from vertex~$y$ to vertex~$z$, then we must have $D(0,y)\le j$ and~$D(z,x)\le k-j$. Conditioning on~$j$ and $(y,z)$ thus yields
\begin{equation}
\BbbP\bigl(D(0,x)\le k\bigr)\le\sum_{j=1}^k \!\!\sum_{\begin{subarray}{c}
y,z\in\Z^d\\|y-z|\ge|x|/k
\end{subarray}}\!\!
\BbbP\bigl(D(0,y)\le j\bigr)\,p_{yz}\,\BbbP\bigl(D(z,x)\le k-j\bigr).
\end{equation}
Under the assumption that $|x|/k\ge R$ we can bound $p_{yz}\le(|x|/k)^{-s'}$. Dropping the condition on $|y-z|$ we can now sum over~$y$ and~$z$ to get the right-hand side of \eqref{PB-bd}.
\end{proofsect}

\begin{lemma}
\label{lemma-step2}
There is a constant~$a=a(d,s')$ such that, given $j\ge1$, if there is a $K\ge Rj$ such that $\BbbP(D(0,x)\le j)\le [K/|x|]^{s'}$ for all~$x\in\Z^d$ with $|x|/j\ge R$, then $\E|B_j|\le a K^d$.
\end{lemma}

\begin{proofsect}{Proof}
Note that $|x|> K$ implies $|x|/j\ge R$. Thus
\begin{equation}
\E|B_j|=\sum_{x\in\Z^d}\BbbP\bigl(D(0,x)\le j\bigr)
\le\sum_{x\colon|x|\le K}1+\sum_{x\colon|x|>K}\biggl(\frac{K}{|x|}\biggr)^{s'}.
\end{equation}
It is easy to check that the first term is bounded by a constant $a_1=a_1(d)$ times $K^d$, while the sum over $|x|^{-s'}$ over $|x|>K$ is at most a constant $a_2=a_2(d,s')$ times $K^{d-s'}$. Putting all terms together, the desired claim follows.
\end{proofsect}

In addition to the above lemmas, the proof will require one unpleasant calculation that we formalize as follows:

\begin{lemma}
\label{lemma-step3}
For each $p>\frac{s'+1}{2d-s'}$ and each $c_0>0$ there is $C=C(p,c_0)\in(0,\infty)$ such that for each~$c\ge c_0$ the quantity
\begin{equation}
K(n):=\frac1C(n+1)^{-p}\texte^{\,c\, n^{1/\Delta'}},
\end{equation}
where $\Delta'$ is as above, obeys
\begin{equation}
\label{lhs}
\sum_{j=0}^n K(j)^dK(n-j)^d\le n^{-s'}K(n)^{s'},\qquad n\ge1.
\end{equation}
\end{lemma}

\begin{proofsect}{Proof}
Consider the function $\varphi(x):=x^{1/\Delta'}+(1-x)^{1/\Delta'}$ and note that the exponentials in $K(j)^dK(n-j)^d$ combine into $\exp\{c n^{1/\Delta'} \varphi(j/n)d\}$. Note also that $\varphi$ is maximized at $x:=\ffrac12$ at where it equals $2^{1-1/\Delta'}=s'/d$. Let
\begin{equation}
\delta:=s'-d\max_{0\le x\le1/4}\varphi(x)
\end{equation}
and observe that $\delta>0$. Splitting the sum over~$j$ into the part when $|j-n/2|\le n/4$ or not, and using the symmetry $j\leftrightarrow n-j$ we thus get
\begin{equation}
\begin{aligned}
\sum_{j=0}^n K(j)^dK(n-j)^d
&\le 2\sum_{j\le n/4}K(j)^dK(n-j)^d+\sum_{j\colon|n/2-j|\le n/4}K(j)^dK(n-j)^d
\\
&\le 2\sum_{j\le n/4}C^{-2d}\frac{\texte^{\,c\,n^{1/\Delta'}(s'-\delta)}}{(j+1)^{pd}(n-j+1)^{pd}}
\\
&\qquad\qquad\qquad
+\sum_{j\colon|n/2-j|\le n/4}C^{-2d}\frac{\texte^{\,c\,n^{1/\Delta'}s'}}{(j+1)^{pd}(n-j+1)^{pd}}.
\end{aligned}
\end{equation}
Using  that $j+1\ge (n+1)/8$ and $n-j+1\ge (n+1)/8$ for all integers $j$ such that $|j-n/2|\le n/4$, we now get
\begin{equation}
\begin{aligned}
\text{LHS of \eqref{lhs}}&
\le 2(n+1)C^{-2d}\texte^{\,c\, n^{1/\Delta'}(s'-\delta)}+C^{-2d}8^{2pd}(n+1)^{1-2pd}\texte^{\,c\,n^{1/\Delta'}s'}
\\
&\le
h(n)\,n^{-s'}\Bigl[\frac1C(n+1)^{-p}\texte^{\,c\,n^{1/\Delta'}}\Bigr]^{s'}
\end{aligned}
\end{equation}
where
\begin{equation}
h(n):=C^{s'-2d}\bigl(8^{2pd}(n+1)^{1-2pd}+2(n+1)\texte^{\,-c\delta\,n^{1/\Delta'}}\bigr)(n+1)^{s'+ps'}.
\end{equation}
It is easy to check that $1-2pd+s'+ps'<0$ under the assumed condition on~$p$ and so the term multiplying $C^{s'-2d}$ is bounded uniformly in~$n$ for all $c\ge c_0$. Choosing~$C$ sufficiently small, we get $h(n)\le1$ for all~$c\ge c_0$. This proves the claim.
\end{proofsect}

\begin{proofsect}{Proof of Theorem~\ref{thm3.1}}
Let $p>\frac{s'+1}{2d-s'}$, set $q:=\frac2{2d-s'}$ and let~$a=a(d,s')$ be as in Lemma~\ref{lemma-step2}. Pick~$c_0>0$ and let $C(p,c_0)$ be as in Lemma~\ref{lemma-step3}. Finally, pick $c\ge c_0$ so large that
\begin{equation}
\label{Knbd}
K(n):=\frac1{C(p,c_0)}(n+1)^{-p}\texte^{\,c\,n^{1/\Delta'}}\ge a^q Rn,\qquad n\ge1.
\end{equation}
We will show by induction that, for each~$n\ge1$,
\begin{equation}
\label{my-bd1}
\BbbP\bigl(D(0,x)\le n\bigr)\le \biggl(\frac{a^{-q}K(n)}{|x|}\biggr)^{s'}.
\end{equation}
Notice that this is trivially true for $|x|< a^{-q}K(n)$ and so we may thus always suppose that $|x|\ge a^{-q}K(n)$ which by \eqref{Knbd} implies $|x|\ge Rn$. Also, we may assume that $x\ne0$; otherwise the right-hand side is infinity.

To start the induction we note that \eqref{my-bd1} holds for $n=1$ as, for $x$ away from the origin, $\BbbP(D(0,x)\le1)=p_{0x}$ which is less than the right-hand side by \eqref{2.3a} and the bound $a^{-q}K(1)\ge R\ge1$. So let us now suppose \eqref{my-bd1} holds for all~$n\le m\in\{1,2,\dots\}$ and let us prove it for~$n:=m+1$. Notice that as we may assume $|x|\ge R(m+1)\ge Rj$ for $j=0,\dots,m+1$, Lemma~\ref{lemma-step2} can be used for $\E|B_j|$ with $K:=a^{-q}K(j)$ for all $j=1,\dots,m+1$. By Lemma~\ref{lemma-step1} and Lemma~\ref{lemma-step2} we thus get
\begin{equation}
\BbbP\bigl(D(0,x)\le m+1\bigr)\le\Bigl(\frac{|x|}{m+1}\Bigr)^{-s'}a^{2-2dq}\sum_{j=0}^{m+1} K(j)^d K(m+1-j)^d.
\end{equation}
Invoking Lemma~\ref{lemma-step3}, the sum can be further bounded with the result
\begin{equation}
\BbbP\bigl(D(0,x)\le m+1\bigr)\le a^{2-2dq}\biggl(\frac{K(m+1)}{|x|}\biggr)^{s'}.
\end{equation}
Since $2-2dq=-s'q$, we get \eqref{my-bd1} for $n=m+1$. Thus \eqref{my-bd1} holds for all $n\ge1$; choosing $c_1:=a^{-qs'}C(p,c_0)^{-s'}$ and $c_2:=c$ we then get also \eqref{my-bd}.
\end{proofsect}

\begin{remark}
Notice that summing \eqref{eq3.3} over~$L$ along powers of~$2$ yields
\begin{equation}
\label{glb}
\liminf_{|x|\to\infty}\frac{\log D(0,x)}{\log\log|x|}\ge\Delta,
\qquad\BbbP\text{\rm-a.s.}
\end{equation}
i.e., a lower bound on the growth of the graph distance proved along a far more elegant argument than the original proof in~\cite{Biskup}. 
\end{remark}

\section{Diameter control}
\label{sec4}\noindent
We now pass to the proof of Theorem~\ref{thm1.1}. As remarked earlier, the lower bound in \eqref{lim} is an easy consequence of the asymptotic~\eqref{Dxy}. 

\begin{proofsect}{Proof of Theorem~\ref{thm1.1}, lower bound}
Recall that~$D_L(x,y)$ is the graph distance between~$x$ and~$y$ measured on~$\G_L$ and let~$D(x,y)$ be the distance measured on the full long-range percolation graph on~$\Z^d$. Then we have
\begin{equation}
D_L\ge D_L(x,y)\ge D(x,y),\qquad x,y\in\G_L.
\end{equation}
Now by \eqref{glb} (or~\cite[Theorem~1.1]{Biskup}), for every~$\epsilon>0$, we have $D(0,x)\ge (\log L)^{\Delta-\epsilon}$ once~$L$ is sufficiently large and~$|x|\approx L$. The lower bound in \eqref{lim} follows.
\end{proofsect}

The key is thus to prove the corresponding upper bound. A natural idea is to follow the strategy of~\cite{Biskup} which is based on the following observation: Let~$x$ and~$y$ be two vertices and let~$L:=|x-y|$. Abbreviate
\begin{equation}
B_\ell(x):=x+[-\ell,\ell]^d\cap\Z^d.
\end{equation}
The probability that~$B_\ell(x)$ and $B_\ell(y)$ are directly connected by an edge is then
\begin{equation}
\label{edge-probab}
1-\exp\{-\ell^{2d}L^{-s+o(1)}\}.
\end{equation}
 Thus, choosing $\ell:=L^\gamma$, the aforementioned edge with be present with very high probability as long as $\gamma\gtrapprox\ffrac s{2d}$.
 
Denoting $z_{00}:=x$ and $z_{11}:=y$, and letting $z_{01}$ and~$z_{10}$ be the endpoints of the primary edge $(z_{01},z_{10})$, we can now find two secondary edges of length order $L^{\gamma}$ spanning the ``gaps'' $(z_{00},z_{01})$ and $(z_{10},z_{11})$ to within $L^{\gamma^2}$ from the respective endpoints. Next we identify 4 tertiary edges of length~$L^{\gamma^2}$ that leave behind 8 ``gaps'' of length $L^{\gamma^3}$, etc. In~\cite{Biskup} it was shown that this edge-identification procedure can be iterated $k$-times with $k\approx(\log\log L)/\log(\ffrac1\gamma)$ until we are down to $2^k$ ``gaps'' of size $(\log L)^{o(1)}$. Using the underlying $\Z^d$-lattice structure, we readily extract a path from~$x$ to~$y$ of length $2^k(\log L)^{o(1)}$. Taking~$\gamma\uparrow\ffrac s{2d}$ along with $L\to\infty$, this behaves like $(\log L)^{\Delta+o(1)}$.

Unfortunately, as remarked earlier, this is \emph{not} going to work for controlling the length of paths between \emph{all} pairs of vertices~$x,y\in\Lambda_L$. The reason is that, to construct a path between two fixed points we only need to ensure the presence of $(\log L)^{\Delta+o(1)}$ edges but, to do this uniformly for all pairs of points in~$\Lambda_L$ we would need to control order~$L^{d+o(1)}$ of them. This entropy cannot be beaten since there are blocks of side $(\log L)^{1/d+o(1)}$ with no incident long edges at all. We will thus have to deal with the cases where the requisite connections fail to occur by methods of nearest-neighbor percolation. 

For better understanding of what is to follow, it is actually worth noting that the above strategy is in the least capable of proving a polylogarithmic bound on~$D_L$. Indeed, the identification of successive levels can possibly fail at stage~$k$ only if \emph{somewhere} in $\Lambda_L$ there are two vertices at distance $L^{\gamma^k}$ whose neighborhoods of size $L^{\gamma^{k+1}}$ are not connected by an edge in $\G_L$. By \eqref{edge-probab}, this has probability bounded by
\begin{equation}
L^{2d}\exp\bigl\{-L^{\gamma^k(2d\gamma-s+o(1))}\bigr\}.
\end{equation}
Thus, as long as $L^{\gamma^k}\ge (\log L)^{\theta}$, where $1/\theta<2d\gamma-s$, this will not happen with probability tending to one. Halting the procedure at this step shows that
\begin{equation}
D_L\le (\log L)^{\Delta+\theta+o(1)}.
\end{equation}
Further improvement can be achieved if from this point on we make the successive scales related not by exponent $\gamma$, but by an exponent $\zeta$ which is taken close to one. The procedure can then be made to work up to the point when the gaps are at most of size $(\log L)^{1/(2d-s)+o(1)}$. This is still way too large to infer the desired bound on~$D_L$, but now the gaps are themselves much smaller than $(\log L)^\Delta$ --- see Lemma~\ref{lemma1a} below. It thus remains to show that such bad regions will not come close to one another. This is the content of Proposition~\ref{prop2b} below and this is where methods of nearest-neighbor percolation need to be brought into play.

\smallskip
Having outlined the general strategy, we now turn to the details. Fix an $\epsilon>0$. We will need numbers, $s'$, $\gamma$, $\zeta$ and~$\eta$ subject to the restrictions:
\begin{equation}
\label{2.4}
s<s'<2d\quad\text{and}\quad
\frac{s'}{2d}<\gamma<1,
\end{equation}
\begin{equation}
\frac{\log 2}{\log(1/\gamma)}<\Delta+\epsilon,
\end{equation}
\begin{equation}
\label{zetabd}
\gamma<\zeta<1
\quad\text{and}\quad
\Delta>\frac1{2d\zeta-s'}
\end{equation}
and
\begin{equation}
\label{etadef}
\Delta>\eta>\frac1{2d\zeta-s'}.
\end{equation}
To see that such choices can be made, we note the following relation:

\begin{lemma}
\label{lemma1a}
Let~$s\in(d,2d)$ and let~$\Delta$ be as in in \eqref{Delta}. Then~$\Delta>\frac1{2d-s}$.
\end{lemma}

\begin{proofsect}{Proof}
It suffices to note that $s\mapsto(2d-s)\Delta$ is strictly increasing on $(d,2d)$ and equal to~$d$ (which is at least one) at~$s=d$. For this, write~$(2d-s) \Delta=(2d\log2)/f(\ffrac s{2d})$ with $f(x):=\frac1{1-x}\log(\ffrac1x)$. A computation shows that~$f'(x)<0$ for $0<x<1$.
\end{proofsect}

Using $\vee$ to denote the maximum (and $\wedge$ the minimum), with the above~$s'$, $\gamma$, $\eta$ and~$\zeta$ fixed, we also choose a quantity $\theta$ such that
\begin{equation}
\label{thetadef}
\theta>\frac1{2d\gamma-s'}\vee\eta
\end{equation}
and define
\begin{equation}
k_0:=\max\bigl\{k\ge1\colon \lfloor L^{\gamma^k}\rfloor>(\log L)^\theta\bigr\}.
\end{equation}
For any (large) positive integer~$L$ we now define a family of scales~$(L_k)$ as follows:
\begin{equation}
L_k:=\begin{cases}
\lfloor L^{\gamma^k}\rfloor,\qquad&\text{if }k\le k_0,
\\*[2mm]
\lfloor L^{\gamma^{k_0}\zeta^{k-k_0}}\rfloor,\qquad&\text{otherwise}.
\end{cases}
\end{equation}
Thus, for~$k\le k_0$ the subsequent scales are related by exponent~$\gamma$, while beyond~$k_0$ the corresponding exponent is ``only''~$\zeta$. In particular, the subsequent scales for~$k>k_0$ are far closer than for~$k\le k_0$.
For later purposes we will need to introduce other two distinguished values:
\begin{equation}
k_1:=\max\bigl\{k\ge1\colon L_k>(\log L)^{\eta}\bigr\}
\end{equation}
and
\begin{equation}
k_2:=\min\bigl\{k\ge1\colon L_k<(\log L)^{\epsilon}\bigr\}.
\end{equation} 
We will only need to consider the scales $L_k$ up to $k=k_2$.
A forthcoming definition (of good blocks) will also depend on a~$\delta>0$ that is picked so small that
\begin{equation}
\label{delta-epsilon}
(1-\delta)^{k_2-k_1}>\ffrac12.
\end{equation}
This is possible since $k_2-k_1$ is bounded by a constant times $\log(\ffrac\eta\epsilon)/\log(1/\zeta)$.

\smallskip
Now consider the cubic box~$\Lambda_L$. We wish to partition~$\Lambda_L$ into blocks of scale~$L_1$ which in turn should be partitioned into blocks of scale~$L_2$, etc. Unfortunately, the subsequent scales may not be divisible by one another and so we will have to work with rectangular boxes of uneven dimensions. For an integer $\ell\ge1$, we call an \emph{$\ell$-block} any translate of
\begin{equation}
\bigl\{(n_1,\dots,n_d)\colon 1\le n_i\le \ell_i,\,i=1,\dots,d\bigr\},
\end{equation}
where $\ell_1,\dots,\ell_d$ are numbers such that $\ell/2\le\ell_i\le\ell$ for all~$i=1,\dots,d$. We note:

\begin{lemma}
\label{lemma1b}
Let~$0<\ell'<\ell$ be integers. Then any~$\ell$-block can be partitioned into~$\ell'$-blocks.
\end{lemma}

\begin{proofsect}{Proof}
Since the partitioning can be done independently in each lattice direction, we may assume $d=1$. Without loss of generality, let~$\ell$ be the actual size of the larger block. Let~$n$ be the unique integer such that~$n\ell'<\ell\le(n+1)\ell'$. If~$(n+\frac12)\ell'<\ell$, then $\ell-n\ell'\in[\ell'/2,\ell']$ and we may decompose the~$\ell$-block into $n$ blocks of side~$\ell'$ and one block of side~$\ell-n\ell'$. If instead $(n+\frac12)\ell'\ge\ell$, then we use only~$n-1$ blocks of side~$\ell'$ and two blocks of about the same size between~$\ell'/2$ and~$\ell'$ whose combined side-length is~$\ell-(n-1)\ell'$ --- a number between $\ell'$ and $\frac32\ell'$.
\end{proofsect}

Given~$L$, we will now choose a partitioning of~$\Lambda_L$ into $L_1$-blocks, a partitioning of these $L_1$-blocks into $L_2$-blocks, etc, for all scales $L_k$ with $k\le k_2$. The hierarchical decomposition will be fixed for the remainder of the argument.

\smallskip
Next we designate good and bad blocks as follows:

\begin{definition}[Good/bad blocks]
\label{def-good}
For the above hierarchical decomposition of~$\Lambda_L$, define good blocks as follows:
\settowidth{\leftmargini}{(111a)}
\begin{enumerate}
\item[(1)]
Any $L_{k_2}$-block is good.
\end{enumerate}
If~$k< k_2$, an $L_k$-block is said to be good if
\settowidth{\leftmargini}{(1111a)}
\begin{enumerate}
\item[(2a)]
at least $1-\delta$ fraction of the $L_{k+1}$-blocks contained therein are good, and
\item[(2b)]
any two distinct good $L_{k+1}$-subblocks are linked by an edge from~$\G_L$ whose endpoints are contained only in good $L_{k'}$-blocks, for all~$k'=k+1,\dots,k_2$.
\end{enumerate}
An $L_k$-block is \emph{bad} if it is not good.
\end{definition}

Let $\scrB_k$ be the union of all bad $L_k$-blocks and let
\begin{equation}
\scrB:=\bigcup_{k=k_1}^{k_2}\scrB_k.
\end{equation}
We will refer to vertices in~$\Lambda_L\setminus\scrB$ as \emph{good} and those in~$\scrB$ as~\emph{bad}.
The nearest-neighbor structure on~$\Z^d$ induces a decomposition of~$\scrB$ into connected components; let~$C(x)$ denote the connected component of~$\scrB$ that contains~$x$ and define
\begin{equation}
T_L(x):=\diam C(x) \quad\text{and}\quad T_L:=\max_{x\in\Lambda_L}T_L(x).
\end{equation}
Next, consider the restriction~$\G_L'$ of~$\G_L$ to vertex set~$\Lambda_L\setminus\scrB$ and let $D'_{L}(x,y)$ be the graph-theoretical distance between~$x$ and~$y$ as measured on~$\G_L'$. Define
\begin{equation}
D_L':=\max_{x,y\in\Lambda_L\setminus\scrB}D'_{L}(x,y).
\end{equation}
Notice that $T_L$ and~$D_L'$ depend on the choices of $s',\gamma,\eta,\delta$ and~$\epsilon$.
Our proof of the upper bound in \eqref{lim} is now reduced to the following propositions:

\begin{proposition}
\label{prop2a}
For~$\epsilon$ as above, 
\begin{equation}
\lim_{L\to\infty}\BbbP\bigl(\,D_L'\le(\log L)^{\Delta+2\epsilon}\bigr)=1.
\end{equation}

\end{proposition}

\begin{proposition}
\label{prop2b}
For~$\epsilon$ as above, 
\begin{equation}
\label{2.7a}
\lim_{L\to\infty}\BbbP\bigl(\,T_L\le(\log L)^{\Delta+\epsilon}\bigr)=1.
\end{equation}
\end{proposition}

These are proved in the next section. Subject to these propositions, we are now ready to establish the main result of this work:

\begin{proofsect}{Proof of Theorem~\ref{thm1.1}, upper bound}
Pick $x,y\in\Lambda_L$. If~$x$ is contained in a bad block then, within $\Z^d$-distance~$T_L(x)$, there is a vertex~$x'\in\Lambda_L\setminus\scrB$, and similarly we find a vertex~$y'\in\Lambda_L\setminus\scrB$ within distance~$T_L(y)$ of~$y$. By concatenating the shortest path between~$x'$ and~$y'$ on~$\G_L'$ with shortest paths connecting~$x$ to~$x'$ and $y$ to~$y'$ on~$\Z^d$, we have
\begin{equation}
D_L(x,y)\le T_L(x)+T_L(y)+D'_{L}(x',y').
\end{equation}
Therefore,
\begin{equation}
D_L\le 2T_L+D_L'.
\end{equation}
By Propositions~\ref{prop2a} and~\ref{prop2b}, the right hand side is bounded by $3(\log L)^{\Delta+2\epsilon}$ with probability tending to one. As $\epsilon$ was arbitrary positive, the claim follows. 
\end{proofsect}

\section{Taming the bad blocks}
\label{sec5}\noindent
To finish the proof of Theorem~\ref{thm1.1} we have to provide proofs of Propositions~\ref{prop2a} and~\ref{prop2b}. We begin by a lemma. Recall that a vertex $x\in\Lambda_L$ is good if it is contained only in good $L_k$-blocks, for all $k=k_1,\dots,k_2$. Then we have:

\begin{lemma}
\label{lemma-good}
For each~$k=k_1,\dots,k_2$, each good $L_k$-block contains at least half of good vertices. In addition, if $L_k\ge32dR$, at least quarter of the vertices in the $L_k$-block with distance at least~$R$ from the boundary are good.
\end{lemma}

\begin{proofsect}{Proof}
We claim that, in fact, at least $(1-\delta)^{k_2-k}$ fraction of all vertices in a good $L_k$-block are good. This is obviously true for $k=k_2$, as all $L_{k_2}$-blocks are good by  Definition~\ref{def-good}(1). For $k=k_1,\dots,k_2-1$ this is proved by induction using Definition~\ref{def-good}(2a). The first part of the claim now follows by invoking the bound \eqref{delta-epsilon}. 

To get the second part, we note that for each lattice direction, at least $4R/L_k$-fraction of all vertices are closer than~$R$ to the sides of the block in this direction. Thus less than $8dR/L_k\le\ffrac14$ of all vertices in the $L_k$-block are at least~$R$-away from any side. If half of all vertices in the $L_k$-block are good, at least quarter of all vertices at least~$R$-away from the boundary must be good.
\end{proofsect}

For the probability estimates that are to follow, it will be useful to note that by \eqref{plim}, for each~$s'\in(s,2d)$ there is a number $R=R(s')<\infty$ such that
\begin{equation}
\label{2.3}
p_{xy}\ge 1-\texte^{-|x-y|^{-s'}},\qquad |x-y|\ge R.
\end{equation}
Note that this is different from \eqref{2.3a}, where we cared for an upper bound on $p_{xy}$.

\begin{proposition}
\label{prop-bad-bd}
Given an $L_k$-block, let $\EE_k$ be the event that this block is good. There are constants~$c_1,c_2\in(0,\infty)$ such that, whenever $L$ is so large that $L_{k_2}\ge 8d R$, we have
\begin{equation}
\label{EE-bound}
\BbbP(\EE_k^\cc)\le  c_1\texte^{-c_2L_k^{2d\zeta-s'}},\qquad k=k_1,\dots,k_2.
\end{equation}
\end{proposition}

\begin{remark}
The above estimate shows why we need to make the subsequent scales $L_k$ related by exponent $\zeta$ --- which can be taken arbitrarily close to one --- and not $\gamma$ (as is done for scales for $k<k_1$). Indeed, the bound \eqref{EE-bound} permits the existence of bad $L_k$ blocks already when $L_k^{2d\zeta-s'}=(\log L)^{1+o(1)}$. When $\zeta$ obeys \eqref{zetabd}, this rules out existence of bad $L_k$ blocks with $L_k\approx(\log L)^\Delta$, but if we worked with $\zeta=\gamma$ this (and consequently, Proposition~\ref{prop2b}) would fail once $\gamma$ is not sufficiently close to one. Another instance where the difference between~$\zeta$ and~$\gamma$ shows up is the derivation \twoeqref{G-F-bd}{cup-G-F-bd}.
\end{remark}

Proposition~\ref{prop-bad-bd} will be established by proving a recursive estimate on the probability that an $L_k$-block is bad:

\begin{lemma}
Let $a_k$ be the maximum value of $\BbbP(\EE_k^\cc)$ over all $L_k$-blocks. There are constants $c_3,c_4,c_5\in(0,\infty)$ such that when $L_{k_2}\ge32dR$, the sequence $(a_k)$ obeys the recursive bound
\begin{equation}
\label{recur-bd}
a_k\le(2a_{k+1})^{c_3 L_k^{2d(1-\zeta)}}+c_4 L_k^{2d}\texte^{-c_5 L_k^{2d\zeta-s'}},\qquad k=k_1,\dots,k_2-1,
\end{equation}
with terminal condition
\begin{equation}
\label{ak1-nula}
a_{k_2}:=0.
\end{equation}
\end{lemma}

\begin{proofsect}{Proof}
Pick an $L_k$-block and let $\AA_k$ be the event that at least $1-\delta$ fraction of all $L_{k+1}$-blocks therein is good. Then we can bound $\BbbP(\EE_k^\cc)$ by
\begin{equation}
\BbbP(\EE_k^\cc)\le\BbbP(\AA_k^\cc)+\BbbP(\EE_k^\cc|\AA_k).
\end{equation}
We will now prove that the probabilities on the right-hand side are bounded, respectively, by the two terms in \eqref{recur-bd}.

Let $n_k$ denote the number of $L_{k+1}$-blocks in the given $L_k$-block. By induction assumption $\BbbP(\EE_{k+1}^\cc)$ is bounded by $a_{k+1}$ for each $L_{k+1}$-block and the events~$\EE_{k+1}$ for distinct blocks are independent. It follows that the number of bad $L_{k+1}$-blocks in the given $L_k$-block is stochastically dominated by a binomial random variable with parameters $n_k$ and~$a_{k+1}$. In particular, $\BbbP(\AA_k^\cc)$ is less than the probability that this random variable is at most $\delta n_k$. The exponential Chebyshev bound now gives
\begin{equation}
\BbbP(\AA_k^\cc)\le\texte^{-\lambda\delta n_k}(1-a_{k+1}+\texte^\lambda a_{k+1})^{n_k},\qquad \lambda\ge0.
\end{equation}
Choosing $\texte^{-\lambda}:=a_{k+1}$ and noting that $n_k$ is at least a constant times $(L_k/L_{k+1})^d$ which is bounded by a constant times $L_k^{2d(1-\zeta)}$ then yields
\begin{equation}
\BbbP(\AA_k^\cc)\le(2a_{k+1}^\delta)^{n_k}\le(2a_{k+1})^{c_3 L_k^{2d(1-\zeta)}}.
\end{equation}
This proves the first term on the right-hand side of \eqref{recur-bd}.

To get the second term, we note that $\AA_k$ is determined only by the edges with both endpoints in the same $L_{k+1}$-block. Thus, conditioning $\EE_k^\cc$ on $\AA_k$ means that the set of good vertices in one of the good $L_{k+1}$-blocks contained therein is not joined by an edge to the set of good vertices in another such good $L_{k+1}$-block. As, by Lemma~\ref{lemma-good}, at least a quarter of all vertices in each such good $L_{k+1}$-block that are $R$-away from its boundary are good, the use of \eqref{2.3} permissible and so we have
\begin{equation}
\BbbP(\EE_k^\cc|\AA_k)\le \binom{n_k}2\exp\Bigl\{-\frac{\bigl(\frac14 (L_{k+1}/2)^d\bigr)^2}{(dL_k)^{-s'}}\Bigr\}.
\end{equation}
Here the binomial coefficient counts the number of pairs of $L_{k+1}$-blocks, $(L_{k+1}/2)^d$ is a lower bound on the size of any $L_{k+1}$-block, the factor $\ffrac14$ accounts for the number of good vertices in such $L_{k+1}$-block that are at least distance~$R$ from the boundary and $dL_k$ is the maximum of $|x-y|$ for any pair of such good vertices in the $L_k$-block. Using that $L_{k+1}\ge c^{-1}L_k^\zeta$ and $n_k\le c L_k^d$ for some constant~$c\in(0,\infty)$, the second term on the right-hand side of \eqref{recur-bd} is proved too.
\end{proofsect}

\begin{proofsect}{Proof of Proposition~\ref{prop-bad-bd}}
We have to show how to get \eqref{EE-bound} from \eqref{recur-bd}. Here we we invoke the inequality
\begin{equation}
L_{k+1}^{2d\zeta-s'}L_k^{2d(1-\zeta)}\ge c L_k^{2d\zeta-s'},
\end{equation}
valid for some constant $c\in(0,\infty)$ for all~$k$, to check that an upper bound of the form 
\begin{equation}
a_k\le c_6^{k_2-k}\texte^{-c_2 L_k^{2d\zeta-s'}}
\end{equation}
propagates under this recursion once $c_2$ and $c_6$ are taken sufficiently small but positive. As this bound holds for $k=k_2$ by \eqref{ak1-nula}, it holds for all $k=k_1,\dots,k_2$. Noting that $k_2-k\le k_2-k_1$ is bounded, we have \eqref{EE-bound} with $c_1:=c_6^{k_2-k_1}\vee1$.
\end{proofsect}

An immediate consequence of the bound in Proposition~\ref{prop-bad-bd} is:

\begin{corollary}
\label{cor-F}
Let $\FF_L$ be the event that all $L_{k_1}$-blocks are good. Then $\lim_{L\to\infty}\BbbP(\FF_L)=1$.
\end{corollary}

\begin{proofsect}{Proof}
As the number of $L_{k_1}$-blocks is at most $cL^d$ for some $c<\infty$, a union bound yields
\begin{equation}
\BbbP(\FF_L^\cc)\le cL^d\, c_1\texte^{-c_2 L_{k_1}^{2d\zeta-s'}}.
\end{equation}
By \eqref{etadef} and the definition of $k_1$, the exponent is much larger than $\log L$.
\end{proofsect}

\begin{lemma}
Let $\GG_k$ be the event that in every $L_k$-block, any two distinct $L_{k+1}$-blocks are connected by an edge in $\G_L$ with both endpoints at good vertices. Then
\begin{equation}
\lim_{L\to\infty}\BbbP\biggl(\,\bigcap_{k=0}^{k_1-1}\GG_k\biggr)=1.
\end{equation}
\end{lemma}

\begin{proofsect}{Proof}
Consider the event $\GG_k^\cc\cap\FF_L$ for $0\le k\le k_1$. On this event, Lemma~\ref{lemma-good} ensures that at least half of all vertices in $L_{k+1}$-blocks are good. Moreover, for~$L$ sufficiently large, at least a quarter of the vertices in $L_{k+1}$-block are further than~$R$ from its boundary. The probability that two such blocks at distance at most $L_k$ are not connected by an edge with both endpoints in~$\scrB^\cc$ is bounded by $\exp\{-c L_{k+1}^{2d}/L_k^{-s'}\}$. In particular,
\begin{equation}
\label{G-F-bd}
\BbbP(\GG_k^\cc\cap\FF_L)\le \tilde c L^{2d}\exp\{-c L_{k+1}^{2d}/L_k^{-s'}\}.
\end{equation}
for some $\tilde c<\infty$. Let
\begin{equation}
\alpha:=\min\bigl\{(2d\zeta-s')\eta,(2d\gamma-s')\eta\bigr\}.
\end{equation}
Examining separately the cases $k\ge k_0$ and $k<k_0$,  we find $L_{k+1}^{2d}/L_k^{-s'}\ge a (\log L)^\alpha$ for some constant $a>0$. Plugging this into \eqref{G-F-bd} we infer
\begin{equation}
\label{cup-G-F-bd}
\BbbP\biggl(\FF_L\cap\bigcup_{k=0}^{k_1-1}\GG_k^\cc\biggr)\le k_1 \tilde cL^{2d}\texte^{-a c(\log L)^\alpha}.
\end{equation}
As $k_1=O(\log\log L)$ and, by \eqref{thetadef} and \eqref{etadef}, $\alpha>1$, the right-hand side tends to zero as $L\to\infty$. The claim now follows by invoking Corollary~\ref{cor-F}.
\end{proofsect}

We are now ready to provide the necessary control of the distance function on $\G_L'$.
The key observation we will need is as follows:

\begin{lemma}
\label{lemma-observe}
Assume $\GG:=\bigcap_{k=0}^{k_1-1}\GG_k$ occurs and let $k(z,z')$ denote the maximal~$k$ such that $z$ and~$z'$ belong to the same $L_k$-block. If $z,z'\not\in\scrB$, then for each $k=k(z,z')+1,\dots,k_2$, the $L_k$-blocks $\Lambda$ and $\Lambda'$ containing $z$ and~$z'$, respectively, are connected by an edge in $\G_L'$.
\end{lemma}

\begin{proofsect}{Proof}
For $k\ge k_1$ this follows by Definition~\ref{def-good} (and the fact that~$z$ and~$z'$ are good), for $k<k_1$ this is implied by the fact that $\GG$ occurs.
\end{proofsect}

\begin{proofsect}{Proof of Proposition~\ref{prop2a}}
Pick $x,y\in\Lambda_L\setminus\scrB$ and assume that the event $\GG$ occurs. We will show that then $\G_L'$ contains a path from~$x$ to~$y$ of length at most $(\log L)^{\Delta+\epsilon}$. 

Let~$\Lambda_0$ and~$\Lambda_1$ be the $L_1$-blocks containing~$x$ and~$y$, respectively. If~$\Lambda_0\ne\Lambda_1$, by Lemma~\ref{lemma-observe}, there is an edge in $\G_L'$ with endpoints~$z_{01}\in\Lambda_0$ and~$z_{10}\in\Lambda_1$. If~$\Lambda_0=\Lambda_1$, we set $z_{01}=z_0$ and~$z_{10}=z_1$. This defines the first level of a hieararchy of vertices and edges. For the next level, denote $z_{00}:=z_0$ and~$z_{11}:=z_1$ and for $\sigma\in\{00,01,10,11\}$ let $\Lambda_\sigma$ be the $L_2$-blocks containing $z_{\sigma}$, respectively. As all of the vertices $z_\sigma$ are good, Lemma~\ref{lemma-observe} ensures the existence of edges $(z_{001},z_{010})$ and $(z_{101},z_{110})$ from~$\G_L'$ between ``good'' vertices $z_{001}$, $z_{010}$, $z_{101}$ and $z_{110}$ in the $L_3$-blocks containing~$z_{000}:=x$, $z_{011}:=z_{01}$, $z_{100}:=z_{10}$ and $z_{111}:=z_{11}$, respectively.

Proceeding by induction along scales until we get to level~$k_2$, we will thus identify a collection of vertices~$(z_\sigma)$, indexed by~$\sigma\in\{0,1\}^{k_2+1}$, such that the following properties hold for each~$k\le k_2$:
\begin{enumerate}
\item[(1)]
$z_\sigma:=x$ if~$\sigma=(0,\dots,0)$ while $z_\sigma:=y$ if~$\sigma=(1,\dots,1)$.
\item[(2)]
For each~$\sigma\in\{0,1\}^{k-1}$, the pair~$(z_{\sigma01},z_{\sigma10})$ is connected by an edge from~$\G_L'$.
\item[(3)]
For each~$\sigma\in\{0,1\}^{k-1}$, the vertices~$z_{\sigma00},z_{\sigma01}$ lie in one of the (good)~$L_{k+1}$-blocks, and similarly for the pair~$(z_{\sigma10},z_{\sigma11})$.
\end{enumerate}
Here $\sigma$ is a hierarchical index and ``$\sigma01$'' denotes a concatenation of the string $\sigma$ with ``$01$.'' The subsequent ``generations'' of the hierarchy are nested via the ``cancellation rules:'' $z_{\sigma00}=z_{\sigma0}$ and $z_{\sigma11}=z_{\sigma1}$. The vertices~$z_\sigma$ are not required to be distinct.

The pair of pairs of vertices $(z_{\sigma00},z_{\sigma01})$ and/or $(z_{\sigma10},z_{\sigma11})$, $\sigma\in\{0,1\}^{k_2-1}$, are contained in the same $L_{k_2}$-block; joining them by shortest paths on~$\Z^d$ we thus construct a path on~$\G_L'$ from~$x$ to~$y$. This path has at most $2^{k_2}-1$ long edges and at most $2^{k_2}$ nearest-neighbor paths on~$\Z^d$ each of which is of length at most $dL_{k_2}$. Hence we get
\begin{equation}
\label{2.15}
D'_{L}(x,y)\le 2^{k_2}-1+2^{k_2}dL_{k_2}.
\end{equation}
Invoking the explicit definitions of $k_2$ and $L_{k_2}$, we find
\begin{equation}
2^{k_2}\ll(\log L)^{\Delta+\epsilon}\quad\text{and}\quad L_{k_2}\le(\log L)^\epsilon.
\end{equation}
The right-hand side of \eqref{2.15} is thus at most $(\log L)^{\Delta+2\epsilon}$, uniformly for all good vertices~$x$ and~$y$. As $\GG$ occurs with probability tending to one, the desired claim follows.
\end{proofsect}

The proof is finished by providing a control of the maximal size of connected components of bad vertices. 

\begin{proofsect}{Proof of Proposition~\ref{prop2b}}
We will derive a uniform bound on the probability $\BbbP(T_L(x)\ge(\log L)^{\Delta+\epsilon})$. Suppose~$x\in\scrB$ and consider the connected component $C(x)$. Then~$C(x)$ is the disjoint union of $L_k$-blocks, $k=k_1,\dots,k_2$, all of which are bad. If we fix one such possible collection, $\bigcup_k\bigcup_{i=1}^{m_k}\Lambda_k^{(i)}$ containing disjoint $L_k$-blocks $\Lambda_k^{(i)}$, $i=1,\dots,m_k$, Proposition~\ref{prop-bad-bd} and the fact that disjoint blocks are independent yields
\begin{equation}
\label{C-Lambda}
\BbbP\biggl(\,C(x)=\bigcup_{k=k_1}^{k_2-1}\bigcup_{i=1}^{m_k}\Lambda_k^{(i)}\biggr)
\le\prod_{k=k_1}^{k_2-1}\biggl\{c_1\texte^{-c_2L_k^{2d\zeta-s'}}\biggr\}^{m_k}
\end{equation}
Now, if $\diam C(x)\ge t$, then $\sum_k m_k L_k\ge t$ and so
\begin{equation}
\sum_{k=k_1}^{k_2-1}m_k L_k^{2d\zeta-s'}\ge t^{1\wedge(2d\zeta-s')}
\end{equation}
As the number of distinct collections of $m:=\sum_k m_k$ blocks that may give rise to $C(x)$ is bounded by $[2d(k_2-k_1)]^{2m}$, we may thus borrow half of the exponent in \eqref{C-Lambda} and use the rest to control the entropy. This yields
\begin{equation}
\BbbP\bigl(T_L(x)\ge t\bigr)
\le\texte^{-\frac12c_2t^{1\wedge(2d\zeta-s')}}
\sum_{m\ge1}\biggl\{4d^2(k_2-k_1)^2\,c_1\texte^{-\frac12c_2L_{k_2}^{2d\zeta-s'}}\biggr\}^m.
\end{equation}
As $k_2-k_1=O(\log\log L)$ while $L_{k_2}\ge(\log L)^{\zeta\epsilon}$, the term in the large braces is small as soon as~$L$ is sufficiently large. Setting $t:=(\log L)^{\Delta+\epsilon}$ and noting that
\begin{equation}
(\Delta+\epsilon)\bigl(1\wedge(2d\zeta-s')\bigr)>1
\end{equation}
by \eqref{zetabd} and/or $\Delta>1$, the probability that $T_L(x)\ge(\log L)^{\Delta+\epsilon}$ is $o(L^{-d})$ uniformly in~$x$. The claim is finished by a standard union bound.
\end{proofsect}

\section*{Acknowledgments}
\noindent
This research was, at various stages, partially supported by the NSF under the grants DMS-0306167 and DMS-0806198/DMS-0949250. I~wish to thank Jianyang Zeng who spotted an error in the initial write-up of this proof and kindly brought it to my attention.

\end{document}